\documentclass{amsart}

\usepackage{latexsym}
\usepackage{epsfig}
\usepackage{amsthm}
\usepackage{amssymb}
\usepackage{amsmath}
\usepackage{amscd}

\newtheorem {Theorem}   {Theorem}

\newtheorem{thm}[subsection]{Theorem}

\newtheorem{prop}[subsection]{Proposition}
\newtheorem{cor}[subsection]{Corollary}
\newtheorem{lemma}[subsection]{Lemma}
\newtheorem{rmk}[subsection]{Remark}

\newcommand{\tMuc}{\widetilde{M}_{\mu,c}}

\theoremstyle{remark}
\newtheorem{Remark}[Theorem]{Remark}

\expandafter\chardef\csname pre amssym.def at\endcsname=\the\catcode`\@
\catcode`\@=11
\def\undefine#1{\let#1\undefined}
\def\newsymbol#1#2#3#4#5{\let\next@\relax
 \ifnum#2=\@ne\let\next@\msafam@\else
 \ifnum#2=\tw@\let\next@\msbfam@\fi\fi
 \mathchardef#1="#3\next@#4#5}
\def\mathhexbox@#1#2#3{\relax
 \ifmmode\mathpalette{}{\m@th\mathchar"#1#2#3}%
 \else\leavevmode\hbox{$\m@th\mathchar"#1#2#3$}\fi}
\def\hexnumber@#1{\ifcase#1 0\or 1\or 2\or 3\or 4\or 5\or 6\or 7\or 8\or
 9\or A\or B\or C\or D\or E\or F\fi}

\font\teneufm=eufm10
\font\seveneufm=eufm7
\font\fiveeufm=eufm5
\newfam\eufmfam
\textfont\eufmfam=\teneufm
\scriptfont\eufmfam=\seveneufm
\scriptscriptfont\eufmfam=\fiveeufm

\catcode`\@=\csname pre amssym.def at\endcsname

\def    \eps    {\epsilon}

\newcommand{\R}{{\bf R}}

\newcommand{\N}{{\bf N}}

\newcommand{\Z}{{\bf Z}}

\newcommand{\CP}{{\bf CP}}

\newcommand{\al}{{\alpha}}

\newcommand{\be}{{\beta}}

\newcommand{\om}{{\omega}}

\newcommand{\ga}{{\gamma}}

\newcommand{\si}{{\sigma}}

\newcommand{\Symp}{{\rm Symp}}

\newcommand{\Tor}{{\rm Tor}}

\newcommand{\Q}{{\bf Q}}

\newcommand{\La}{{\Lambda}}

\newcommand{\id}{{\rm id }}

\newcommand{\MS}{{\medskip}}
\newcommand{\BS}{{\bigskip}}
\newcommand{\NI}{{\noindent}}
\newcommand{\Emb}{{\rm Emb}}

\usepackage{amssymb}
\usepackage[all]{xy}
\usepackage{enumerate}

\def    \R  {{\mathbb R}}
\def    \Z  {{\mathbb Z}}
\def    \N  {{\mathbb N}}
\def    \Q  {{\mathbb Q}}

\def    \CP {{\mathbb C}{\mathbb P}}

\def    \rk {\operatorname{rk}}

\def\qedbox{$\square$}%
\def\qed{\ifmmode\qedbox\else\unskip\ \hglue0mm\hfill
     \qedbox\smallskip\goodbreak\fi}%
\def\QED{\qed\goodbreak\vskip6pt}%

\begin{document}






\title[Symplectic balls in $S^2 \times S^2$ above the critical value]
{The homotopy type of the space of symplectic balls in $S^2 \times S^2$ above the critical value }

\author[S\'{\i}lvia Anjos]{S\'{\i}lvia Anjos}
\author[Fran\c{c}ois Lalonde]{Fran\c{c}ois Lalonde}

\address{ Centro de An\'{a}lise Matem\'{a}tica, Geometria e Sistemas Din\^{a}micos, Departmento de Matem\'{a}tica, Instituto Superior T\'ecnico, Lisbon, Portugal ; 
D\'{e}partement de
math\'{e}matiques et de statistique, Universit\'{e} de Montr\'{e}al.}
\email{sanjos@math.ist.utl.pt ; lalonde@dms.umontreal.ca}

\thanks{The first author is partially supported by FCT/POCTI/FEDER and by project POCTI/2004/MAT/57888.
The second author is partially supported by a Canada Research Chair, NSERC grant OGP 0092913
(Canada) and FQRNT grant ER-1199 (Qu\'ebec).}

\date{}

\begin{abstract}

We compute in this note the full homotopy type of the space of symplectic embeddings of the standard ball $B^4(c) \subset \R^4$ (where $c= \pi r^2$ is the capacity of the standard ball of radius $r$) into the $4$-dimensional rational symplectic manifold $M_{\mu}= (S^2 \times S^2, \mu \om_0 \oplus \om_0)$ where $\om_0$ is the area form on the sphere with total area $1$ and $\mu$ belongs to the interval $(1,2]$. We know, by the work of Lalonde-Pinsonnault \cite{LP}, that this space retracts to the space of symplectic frames of $M_{\mu}$ for any value of $c$ smaller than the critical value $\mu -1$, and that its homotopy type does change when $c$ crosses that value. In this paper, we compute the homotopy type for the case $c \ge \mu-1$ and prove the claim announced in \cite{LP}  that it does not have the type of a finite CW-complex.  
\end{abstract}

\maketitle

\bigskip
\noindent
Subject classification: 53D35, 57R17, 55R20, 57S05.

\section{Introduction and preliminaries} \label{se:intro}

We first recall the  background needed from \cite{LP}. 
Consider the rational symplectic manifold $M_{\mu}= (S^2 \times S^2, \mu \om_0 \oplus \om_0)$ where $\om_0$ is the area form on the sphere with total area $1$ and where $\mu$ belongs to the interval $(1,2]$. The generalization for all values of $\mu$, that is, the cases $\mu \in (n, n+1]$ are treated in the forthcoming paper \cite{ALP}.   
Let
$B^4(c)\subset\R^4$ be the closed standard ball of radius $r$ and capacity $c=\pi r^2$ equipped with the restriction of the symplectic structure $\om_{st}=dx_1\wedge dy_1+dx_2\wedge dy_2$ of $\R^4$. Let $\Emb_{\om}(c, \mu)$ be the space, endowed with the $C^{\infty}$-topology, of all symplectic embeddings of $B^4(c)$ in $M_{\mu}$. Finally, let $\Im\Emb_{\om}(c,\mu)$ be the space of subsets of $M$ that are images of maps belonging to $\Emb_{\om}(c,\mu)$ defined as the topological quotient
\begin{equation} \label{Fibration-1}
\Symp(B^4 (c))\hookrightarrow\Emb_{\om}(c,\mu)
\longrightarrow\Im\Emb_{\om}(c, \mu)
\end{equation}
where $\Symp(B^4(c))$ is the group, endowed with the $C^{\infty}$-topology, of
symplectic diffeomorphisms of the closed ball, with no restrictions on the
behavior on the boundary (thus each such map extends to a symplectic diffeomorphism of a neighborhood of $B^4(c)$ that sends $B^4(c)$ to itself). We may view $\Im\Emb_{\om}(c,\mu)$ as the space of all unparametrized balls of capacity $c$ of $M$.

We know from \cite{LP} that $\Emb_{\om}(c,\mu)$ retracts to the space of symplectic frames of $M=S^2 \times S^2$ for all values $c < \mu - 1$ (hence the space $\Im \Emb_{\om}(c,\mu)$ retracts to $S^2 \times S^2$ in this range of $c$'s). Recall that the Non-squeezing theorem implies that this space is empty for $c \ge 1$. It is shown moreover in \cite{LP} that the homotopy type of $\Im \Emb_{\om}(c,\mu)$ changes only when $c$ crosses the critical value $\mu - 1$ (and that it does change at that value). The goal of the paper is to compute the full homotopy type of $\Im \Emb_{\om}(c,\mu)$ for values $\mu - 1 \le c <  1$ and to prove, in particular, that it does not have the homotopy type of a finite dimensional CW-complex. The main results of this paper are summarised in the following theorem.

\begin{thm} \label{thm:main} Let $0 < \mu - 1 \leq c<1$. Then we have for $\Im \Emb_{\om}(c, \mu)$: 

\smallskip  \NI
1.a) The topological space $\Im \Emb_{\om}(c, \mu)$ is weakly  homotopy equivalent to the total space $E$  of a non-trivial fibration 
\begin{equation}\label{mainfibr}
\Omega \Sigma^2 SO(3) / \Omega S^3 \longrightarrow E \stackrel{\pi}{\longrightarrow}S^2 \times S^2
\end{equation}
where the inclusion of the group $\Omega S^3$ in $\Omega \Sigma^2 SO(3) $ is induced by the map $S^3 \to \Sigma^2 SO(3)$ that corresponds to the generator of the fundamental group of $SO(3)$.  The definition of  $E$ and  $\pi$ is given in \S~\ref{se:fullhomtype}. This fibration  has a continuous section and all its homotopy groups  split.

\smallskip \NI
1.b) The minimal model of $\Im \Emb_{\om}(c, \mu)$ is $(\La (a,b,,e,f,g,h),d)$ with generators in degrees $2,2,3,3,3,4$ and with differential
$$ d e = a^2, \; d f = b^2, \;  d g = d a = d b =0, \; d h = kbg$$
where $k$ is a non zero rational number.
Thus its rational cohomology ring is equal to 
the algebra $$ H^*(\Im \Emb_{\om}(c, \mu);\Q)=\Lambda(a,b,g,gh, \hdots,gh^n, \hdots, bh, \hdots, bh^n, \hdots )/\langle a^2,b^2,bg  \rangle$$
where $n \in \N$. It is therefore not homotopy equivalent to a finite-dimensional CW-complex. 

\smallskip  \NI
{\rm (See \S~\ref{sc:mm} for a geometric description of these generators.)}  
Similarly we have for $\Emb_{\om}(c, \mu)$:

\smallskip \NI
2.a) The topological space $\Emb_{\om}(c, \mu)$ is weakly  homotopy equivalent to the pull-back of the above fibration \eqref{mainfibr}
by the fibration $F_{\omega} \to S^2 \times S^2$ where $F_{\omega}$ is the space of symplectic frames over $M_{\mu} = S^2 \times S^2$.

\smallskip \NI
2.b) A minimal model of $\Emb_{\om}{(c, \mu)}$ is $(\La(\widetilde d_{a,b}, \widetilde e, \widetilde f, \widetilde g, v, \widetilde h), d_0)$ with generators of degrees $2,3,3,3,3,4$ and with differential 
$$d_0 \widetilde d_{a,b}  = d_0 \widetilde f = d_0 \widetilde g = d_0 v =0,  \; d_0 \widetilde h =- k\widetilde d_{a,b} \widetilde g  \mbox{ \ and \ }  d_0 \widetilde e =   \widetilde d_{a,b}^{\; 2}.$$
Its rational cohomology ring is given by 
$$ H^*(\Emb_{\om}{(c, \mu)}; \Q)=\Lambda (\widetilde d_{a,b},\widetilde f,\widetilde g, v, \widetilde g \widetilde h, \hdots,\widetilde g \widetilde h^n, \hdots, s_1, \hdots, s_n, \hdots)/\langle \widetilde d_{a,b}^{\;2},\widetilde d_{a,b}\widetilde g \rangle $$ where $s_n=\widetilde h^{n-1}(\widetilde h \widetilde d_{a,b} + nk \widetilde e \widetilde g)$ and $n \in \N$.  \end{thm} 

\smallskip

See \S~\ref{se:cohomology} for the computations of the cohomology rings with any field coefficients of both $\Im \Emb_{\om}(c, \mu)$ and $\Emb_{\om}(c, \mu)$.

We end this section with a brief description of the methods that we will use.

Denote by $\iota_c$, $c<1$, the standard symplectic embedding of $B^4 (c)$ in
$M_{\mu}$. It is defined as the composition $B^4(c)\hookrightarrow D^2(\mu- \eps)\times D^2(1-\eps)\hookrightarrow S^2(\mu) \times S^2(1) = M_{\mu}$ where the parameters between parentheses represent the areas. Consider the following fibration
\begin{equation}\label{Fibration-2}
\Symp(M_{\mu},B_c)\hookrightarrow\Symp(M_{\mu})
\longrightarrow\Im\Emb_{\om}(c, \mu)
\end{equation}
where the space in the middle is the group of all symplectic diffeomorphisms of
$M_{\mu}$, the one on the left is the subgroup of those that preserve (not
necessarily pointwise) the image $B_c$ of $\iota_c$, and the right hand side arrow assigns to each diffeomorphism $\phi$ the image of the composition $\phi\circ\iota_c$.

Denoting by $\tMuc$ the blow-up of $M_{\mu}$ at the ball $\iota_c$ and by  $\Symp(\tMuc)$ the group of its symplectomorphisms, it is shown in \cite{LP}  that $\Symp(M_{\mu},B_c)$ is naturally homotopy equivalent to $\Symp(\tMuc)$. There is therefore a homotopy fibration
\begin{equation}\label{Fibration-3}
\Symp(\tMuc)\hookrightarrow\Symp(M_{\mu})
\longrightarrow\Im\Emb_{\om}(c, \mu)
\end{equation}
The topology of the middle group has been computed by Abreu \cite{Ab} and by Anjos and Anjos-Granja \cite{An,AG}, while  the topology of the left one was computed in \cite{LP} which also provided the first steps in using that homotopy fibration for the computation of the space $\Im\Emb_{\om}(c, \mu)
$. We will push further their methods to get a complete description of the homotopy type of $\Im\Emb_{\om}(c, \mu)$ and of $\Emb_{\om}(c, \mu)$. We begin with the calculation of the full homotopy type of these two spaces in \S~\ref{se:fullhomtype} and then proceed in \S~\ref{se:rathomtype}  and \S~\ref{se:rathomtypeb} to the computation of their  minimal models.  This ``rational homotopy type'' approach through minimal models will be needed for the cohomology computations of these spaces with any field coefficients in \S~\ref{se:cohomology}.

\bigskip
\NI
{\bf Acknowledgements.}
The first author is grateful to Gustavo Granja for useful conversations, and the second one thanks Octav Cornea for discussions on some aspects of the theory of minimal models. We are also very grateful to the referee of the International Journal of Mathematics for reading the paper so carefully and giving pertinent suggestions, especially for giving a way  to correct the computation of the differential of the element $h$ in the minimal model of  $\Im \Emb_\omega(c,\mu)$ (see Lemma~\ref{diffd}).

\section{The homotopy types of  $\Im\Emb_{\om}(c, \mu)$ and $\Emb_{\om}(c, \mu)$} \label{se:fullhomtype}

In order to find the homotopy type of  $\Im \Emb = \Im\Emb_{\om}(c, \mu)$ we use the knowledge of  the full homotopy  types of $\Symp({\widetilde{M}}_{\mu,c})$ and $\Symp(M_\mu)$ as topological groups. By definition, given topological groups $G$ and $G'$, $G$ has the homotopy type of $G'$ {\em as topological groups} if there are $H$-maps $f:G \to G'$ and $g:G' \to G$ such that $f \circ g$ and $g \circ f$ are homotopic to the identity map through based maps.  Therefore   the product structure is preserved and  the Pontryagin rings of  the two topological groups are isomorphic. 

\medskip

We first recall  from \cite{AG} the background that will be needed.
Let $P$ denote the pushout of the diagram of topological groups 
\begin{equation}
\label{dg:push}
 \xymatrix{ 
SO(3) \ar[r]^-\Delta \ar[d]^i & SO(3) \times SO(3) \\
SO(3) \times S^1 & }
\end{equation}
where $\Delta$ denotes the inclusion of the diagonal and $i$ the inclusion of the first factor.
This is also known as the amalgam (or amalgamated product, or free product with amalgamation) of the groups $SO(3) \times SO(3)$ and $SO(3)\times S^1 $ over the common subgroup $SO(3)$ and it is characterized as the initial topological group admitting compatible homomorphisms from the diagram (\ref{dg:push}).  

Actually,  the pushout construction is defined in any category  (see \cite{Ma}, for instance) but we are mainly interested in the category of topological groups. This construction is characterized by having a universal property that, in our particular case, translates in the following statement. If $G$ is any topological group sitting in the obvious square diagram with continuous homomorphisms making the diagram commutative, then there is a unique continuous homomorphism  $\varphi:P \rightarrow G$ making the following diagram commutative:
$$\xymatrix{
SO(3) \ar[d]^{i} \ar[r]^-{\Delta} & SO(3) \times SO(3)  \ar[d]^{} \ar[ddr]^{} &\\
SO(3) \times S^1 \ar[r]^{} \ar[drr]^{} & P   \ar@{-->}[dr]^{\varphi}  & \\
  &  & G }
$$

In our particular example,  the group structure in the pushout can be described as follows. Consider the set $S$ of finite sequences 
$$x_1x_2\ldots x_n$$ where $x_i \in SO(3) \times SO(3)$ or $x_i \in SO(3) \times S^1$ and  the equivalence relation  $\sim$ generated by
\begin{enumerate}[(i)]
\item $x_1\ldots x_n \sim x_1\ldots \hat{x_i} \ldots x_n $ if $x_i=1$,
\item $x_1 \ldots x_i x_{i+1} \ldots x_n \sim x_1 \ldots (x_i x_{i+1}) \ldots x_n$ when $x_i,x_{i+1}$
both belong either to  $SO(3) \times SO(3)$ or to  $SO(3) \times S^1$.
\end{enumerate}
In particular, if $x_i$ is in the common subgroup $SO(3)$ then $x_1 \ldots (x_{i-1}x_i) x_{i+1} $ $\ldots x_n \sim x_1 \ldots x_{i-1}(x_i x_{i+1}) \ldots x_n$.  $S$ has an associative unital product defined by concatenation and it is easy to check that this
descends to the set $P=S/\sim$ of equivalence classes and that $P$ together with the canonical
maps $SO(3) \times SO(3)\to P$ and $SO(3) \times S^1 \to P$, as a group, is the pushout of diagram (\ref{dg:push}). A good reference for the pushout  of groups is \cite{Se}. 

In  our example, we also need to choose a topology  on the space $P$  so that $P$ together with the maps is the pushout in the category of topological groups. It turns out that it is necessary to work in the category of compactly generated spaces. We say that $X$ is  
{\em compactly generated} if a subset $U \subset X$ is closed iff for every compact Hausdorff space 
$K$ and continuous map $g:K \to X$, $g^{-1}(U)$ is closed. Given an arbitrary topological space $X$ we can refine the topology in the obvious way so that it becomes compactly generated.  Now endow $S$ with the inductive limit of the inclusions $S^n \to S^{n+m}, x_1 \ldots x_n \mapsto  x_1 \ldots x_n \cdot 1 \cdot \ldots \cdot 1$, and $S/ \sim$ with the quotient topology. Finally endow $P$ with the refined topology induced by compact generation. With this topology, $P$ is the push-out of the diagram (\ref{dg:push}) in the category of topological groups. 
A complete description of the pushout of topological groups is given in \cite[\S 2]{AG}.

The homotopy type of  $\Symp(M_\mu)$ was computed by Anjos--Granja in \cite{AG}. To explain their result, recall first that the Hirzebruch surface $W_i$ is given by $$W_i=\{([z_0,z_1],[w_0,w_1,w_2]) \in \mathbb{CP}^1 \times \mathbb{CP}^2 | z_0^iw_1=z_1^iw_0\}.$$ It is well known that the restriction of the projection $\pi_1: \mathbb{CP}^1 \times \mathbb{CP}^2 \to \mathbb{CP}^1 $ to $W_i$ endows $W_i$ with the structure of a K\"ahler $\mathbb{CP}^1$--bundle over $\mathbb{CP}^1$ which is topologically $S^2\times S^2$ if $i$ is even and $\mathbb CP^2 \# \overline{\mathbb CP}^2$ if $i$ is odd. Let's write $B$ for the homology class of the base of the trivial fibration $M= S^2 \times S^2$, i.e. for $[S^2 \times \{pt\}]$, write $F$ for the fiber, i.e. for the class $[\{pt\} \times S^2]$. We will also use the same letters $B,F$ to represent the base and the fiber of $S^2 \times S^2 = B \times F$.
The group $SO(3) \times SO(3)$ can be considered as a subgroup of $\Symp(M_\mu)$ by letting each factor acts on the corresponding factor of $M_{\mu}$ and  the group $SO(3) \times S^1$ is also a subgroup of  $\Symp(M_\mu)$ by 
carrying through the symplectomorphism $W_2 \to M_{\mu}$ the Kahler isometry group   $SO(3) \times S^1$ of the Hirzebruch surface $W_2$ (see \cite{Ab,An,LP}). Here the $S^1$-factor is the rotation in the fibers of $W_2$ -- it is therefore the ``rotation'' in $M_{\mu} = S^2 \times S^2$ of the fibers (of the projection onto the first factor) round symplectic surfaces in class $B+F$ and $B-F$, whereas $SO(3)$ is a lift to $W_2 \to \CP^1= S^2$ of the group $SO(3)$ on the base. This gives two subgroups of $\Symp(M_\mu)$ and it turns out, by a result of Abreu \cite{Ab}, that the first factor of $SO(3) \times S^1$ can be identified with the diagonal in $SO(3) \times SO(3)$. Thus  diagram  (\ref{dg:push}) can be completed with $\Symp(M_\mu)$ in fourth position. This induces a map from $P=(SO(3)\times SO(3))\coprod_{SO(3)} (SO(3) \times S^1)$ to $\Symp(M_\mu)$ by the universal property of pushouts. It was proven in \cite{AG} that 

\begin{thm}[Anjos--Granja] If $0<\mu -1 \leq 1$, then the latter $H$-map 
$$
P=(SO(3)\times SO(3))\coprod_{SO(3)} (SO(3) \times S^1) \to \Symp(M_\mu)
$$
 is a weak homotopy
equivalence.   \end{thm}

Next we compute the homotopy type of $\Symp({\widetilde{M}}_{\mu,c})$ as topological group.
Using the same notations as in \cite{LP},  let the 2-torus $T_i^2$ be the  group of K\"ahler isometries of the blow up $\widetilde{W}_{i,c}$ of the Hirzbruch surface $W_i$ at a standard ball of capacity $c<1$ centered at a point on the zero section of $W_i$ \cite[Prop 4.4 and Prop 4.5]{LP}. Each torus  $T_i^2$ gives rise to an abelian subgroup of   $\Symp({\widetilde{M}}_{\mu,c})$ that we will denote by $\widetilde{T}_i^2$. When $\mu \in  (1,2]$ and $c  \ge \mu-1$, only the tori
$\widetilde{T}_0^2$ and $\widetilde{T}_1^2$ exist. It turns out that the first torus is the product $S^1 \times S^1$, that can be considered as subgroup of the group $SO(3) \times SO(3)$ of diagram (5) (when $\Symp({\widetilde{M}}_{\mu,c})$ is thought of as the subgroup of $\Symp(M_\mu)$ that preserves -- not necessarily pointwise -- the ball of capacity $c$ -- see \cite{LP}). The second torus may be viewed as the subgroup $S^1 \times S^1 \subset SO(3) \times S^1 $ where the second factors are identified and where the first $S^1$-factor is included in $SO(3)$ as the subgroup of the Kahler isometries of $W_2$ that preserves a point on the section at infinity of $W_2$.  This yields the diagram
\begin{equation}
 \xymatrix{ 
S^1 \ar[r]^-\Delta \ar[d]^i &\widetilde{T_0^2} \ar[d] \\
\widetilde{T_1^2} \ar[r] & \Symp({\widetilde{M}}_{\mu,c})  }
\end{equation}

\noindent 
The key point to find the homotopy type of $\Symp({\widetilde{M}}_{\mu,c})$ is to notice that the Pontryagin ring $H_*(\Symp({\widetilde{M}}_{\mu,c}))$ is the pushout of $H_*(\widetilde{T_0^2})$ and $H_*(\widetilde{T_1^2})$ over $H_*(S^1)$, i.e. 

\begin{prop} \label{prop:sq} Let $k$ be a field. The diagram 
\begin{equation}
\label{square}
\xymatrix{
H_*(S^1;k) \ar[r] \ar[d] & H_*(\widetilde{T_0^2};k) \ar[d] \\
H_*(\widetilde{T_1^2};k) \ar[r] & H_*(\Symp({\widetilde{M}}_{\mu,c});k)
}
\end{equation}
is a pushout square in the category of $k$-algebras (there is an entirely analogous description  of the pushout in the category of associative $k$--algebras, see \cite[\S 2]{AG}).
\end{prop}

\proof{}

In order to prove this result we need to compute the Pontryagin ring $H_*(\Symp({\widetilde{M}}_{\mu,c}); k) $ where $k=\Z_p$, with $p$ prime or $k =\Q$. We will show that the generators of the Pontryagin ring with $\Z_p$ coefficients are the same as those with rational coefficients.  The argument to compute these rings follows closely the one of Section 2 in \cite{An}, where the Pontryagin ring $H_*(\Symp(M_\mu);\Z_2)$ is computed. Moreover, we will see that the ring  we want to compute is very similar to the Pontryagin ring  $H_*(\Symp(M_\mu);k)$ when $1< \mu \leq 2$ and $k=\Z_p$, $p>2$ or $k=\Q$, given in Theorem 3.1 in \cite{An}:
$$H_*(\Symp(M_\mu);k) \cong \Lambda(y_3) \otimes k\langle t,x_3 \rangle /R$$
where $\deg y_3=\deg x_3=3$, $\deg t=1$ and $R = \langle t^2, x_3^2\rangle$. In this expression  $k\langle t,x_3 \rangle $ denotes  the free non-commutative algebra over the field $k$ generated by $t,x_3$, and  $\La ( y_3) $ denotes the free graded commutative algebra over the same field generated by $y_3$.
We are going to see that in our case the ring structure is the same except that the generators have all degree 1.  Here we will summarize the whole argument by stating the main steps of that computation. 

Let ${\widetilde{\mathcal{J}}}_{\mu,c}$ be the space of almost complex structures tamed by the symplectic form in $\widetilde{M}_{\mu,c}$. Recall from \cite{LP} that the space $\widetilde{\mathcal{J}}_{\mu,c}$ is a contractible space, and when $0<\mu -1 \leq c <1$  then it is the union of an open, dense, and connected stratum $\widetilde{\mathcal{J}}_0$ with a codimension 2, co-oriented, closed submanifold $\widetilde{\mathcal{J}}_1$. More precisely, each stratum $\widetilde{\mathcal{J}}_i$ is equivalent to the quotient  $\Symp({\widetilde{M}}_{\mu,c})/{\widetilde T}_i $ therefore we have two homotopy fibrations 
$$ {\widetilde T}_i \stackrel{j_i}{\longrightarrow} \Symp({\widetilde{M}}_{\mu,c})\stackrel{p_i}{ \longrightarrow} \widetilde{\mathcal{J}}_i, \ \ \ i=0,1.$$ 
Pinsonnault proves in \cite{Pi} that the inclusion ${\widetilde T}_i \ {\rightarrow} \ \Symp({\widetilde{M}}_{\mu,c})$ induces an injective map in homology for any field coefficients. From the Leray-Hirsch Theorem it follows that the spectral sequence of the fibration collapse at the $E_2$--term, and we obtain the following vector space isomorphisms 
\begin{equation}\label{eq:iso}
H_*( \Symp({\widetilde{M}}_{\mu,c});k) \cong H_*(\widetilde{\mathcal{J}}_i;k) \otimes H_*({\widetilde T}_i;k), \ \ \ i=0,1.
\end{equation}

Let us denote by $x,z$ the generators of the homology of ${\widetilde T}_0$ and by $y,t$ the generators of the homology of the first and second factors of  ${\widetilde T}_1$, respectively. It is easy to see that the generator $y$ commutes with  $z$ since it corresponds to the generator of the diagonal in ${\widetilde T}_0$. It is clear that it also commutes with $t$. 
Recall that for any group $G$ the Samelson product $[x,y]\in \pi_{p+q}(G)$ of elements $x\in \pi_p(G)$ and $y \in \pi_q (G)$ is represented by the commutator 
$$ S^{p+q}=S^p\times S^q/S^p \vee S^q \to G: (u,v) \mapsto x(u)y(v)x(u)^{-1}y(v)^{-1}.$$ The Samelson product in $\pi_*(G)$ is related to the Pontryagin product in $H_*(G, \Z)$ by the formula $$ [x,y]=xy-(-1){|x||y|}yx,$$ where we supressed the Hurewicz homomorphism $\rho:\pi_*(G) \to H_*(G)$ to simplify the expression. One can show that  the commutator $w= [x,t]$, arising from a Samelson product, does not vanish. The argument for the proof of this fact is the same as in section 2.2 of \cite{An}. The next step is to show that the generators of the Pontryagin ring $H_*( \Symp({\widetilde{M}}_{\mu,c});k)$ are in fact $t,x$ and $y$. Again we use the same method used in  \cite{An}, more precisely, in section 2.3.
The main ideia of the proof is to consider the subring $R_* \subset H_*( \Symp({\widetilde{M}}_{\mu,c});k)$ generated  by $t,x$ and $y$ and to assume there is a non-vanishing element in $ H_*( \Symp({\widetilde{M}}_{\mu,c});k) -R_*$. Then, using mainly the isomorphisms \eqref{eq:iso}, we show that this leads to a contradiction. 

Finally we can compute the Pontryagin ring:

\begin{prop}
If $0<\mu -1 \leq c <1$ then
$$ H_*(\Symp({\widetilde{M}}_{\mu,c}); k)= \left (k\langle t,x\rangle /  \widetilde{R} \right ) \otimes \Lambda(y)  $$
where $\widetilde{R}= \langle t^2,x^2 \rangle$, $\deg y=\deg x=\deg t=1$ and $k=\Q$ or $k=\Z_p$ with $p$ prime.
\end{prop}
In order to prove this proposition we note that we already know that the generators of the Pontryagin ring are $t,x$ and $y$. Therefore it is sufficient to show that the only relations between them are the ones in $\widetilde{R}$, the commutativity of $y$ with $x$, and $y^2=0$ coming from the definition of the exterior algebra $\Lambda(y)$. An intermediate step in this proof is to show that the isomorphims \eqref{eq:iso} are in fact given by the Pontryagin product. More precisely, it is not hard to see that we can choose right inverses $s_i:H_*(\widetilde{\mathcal{J}}_i;k) \to  H_*(\Symp({\widetilde{M}}_{\mu,c}); k) $ of ${p_i}_*$, i.e., such that ${p_i}_* \circ s_i= id$. Then one can  show  that the maps
$$ \varphi_i : H_*(\widetilde{\mathcal{J}}_i;k) \otimes  H_*(\widetilde{T}_i;k) \to  H_*(\Symp({\widetilde{M}}_{\mu,c}); k): u\otimes v \mapsto s_i(u).{j_i}_*(v)$$
given by Pontryagin product are isomorphisms. Using these we prove that the only relations in $H_*(\Symp({\widetilde{M}}_{\mu,c}); k)$ are indeed the ones mentioned above. The complete argument can be seen in section 2.4 in \cite{An}. 

 This completes the proof of Proposition~\ref{prop:sq}.   \QED

\begin{Remark}
Notice that in the computation of the Pontryagin ring of $\Symp(M_\mu)$ it is necessary to distinguish between $p=2$ and $p>2$. However, in the computation of the Pontryagin ring of $\Symp({\widetilde{M}}_{\mu,c})$ this distinction is no longer needed since the 2--torus, in contrast to $SO(3)$ does not have 2--torsion.  
\end{Remark}

  Let $\widetilde{P}$ be the pushout of the diagram of topological groups
\begin{equation}
 \xymatrix{ 
S^1 \ar[r]^-\Delta \ar[d]^i &\widetilde{T_0^2}  \\
\widetilde{T_1^2}  &  }
\end{equation}
where $\Delta$ is the inclusion onto the diagonal and $i$ is the inclusion onto the first factor.

We can now state the result for the topological group $\Symp({\widetilde{M}}_{\mu,c})$.
\begin{thm}\label{thm:htMuc}
If $0<\mu -1 \leq c <1$, then the H-map
$$
\widetilde{P}=\widetilde{T_0^2} \coprod_{S^1}\widetilde{T_1^2}  \to \Symp({\widetilde{M}}_{\mu,c})
$$
is a weak homotopy equivalence of topological groups. 
\end{thm}

\proof{}

 First  we need to recall Theorem 3.8 in \cite{AG}. It claims that under the assumption that the homomorphisms involved induce injections  on homology, the following pushout of $k$--algebras 
$$ H_*(B_1;k) \coprod_{H_*(B_0;k)}  H_*(B_2;k)$$ is isomorphic, as an algebra, to the Pontryagin ring $H_*(G;k)$ where $G$ is the amalgamated product of the topological groups $B_1$ and $B_2$ over the common subgroup $B_0$.
It follows from this theorem and proposition \ref{prop:sq} that  the canonical map
$$ \widetilde{P} \rightarrow \Symp({\widetilde{M}}_{\mu,c})$$
is a homology equivalence with any field coefficients. Hence it is an equivalence on integral homology. Since both spaces are $H$--spaces, it follows that  it is in fact a weak equivalence. \QED

Therefore we obtain the following commutative diagram
\begin{equation}
\label{homtype}
\xymatrix{
\widetilde{T_0^2} \coprod_{S^1} \widetilde{T_1^2} \ \ \ar@{^{(}->}[r] \ar[d]^{\tilde{f}} &  (SO(3)\times SO(3))\coprod_{SO(3)} (SO(3) \times S^1) \ar[d]^f \\
\Symp({\widetilde{M}}_{\mu,c})  \ \ \ar@{^{(}->}[r] &  \Symp(M_\mu)
}
\end{equation}
where the vertical arrows are weak homotopy equivalences (as topological groups).

We are now ready to state the main result of this section
\begin{thm}
\label{th:main}
If $0 < \mu -1 \leq c<1$, the topological space $\Im \Emb_{\om}(c, \mu)$ is weakly  homotopy equivalent to the quotient $P/\widetilde P$ of the two pushouts; moreover, this quotient space is the total space of a non-trivial fibration 
$$\Omega \Sigma^2 SO(3)  / \Omega S^3 \longrightarrow \ P/\widetilde P \longrightarrow S^2 \times S^2$$
 where the inclusion of the group $\Omega S^3$ in $\Omega \Sigma^2 SO(3) $ is induced by the map $S^3 \to \Sigma^2 SO(3)$ that corresponds to the generator of the fundamental group of $SO(3)$. This fibration  has a continuous section and splits homotopically, i.e. $$\pi_k(\Im \Emb_{\om}(c, \mu))\simeq \pi_k(\Omega \Sigma^2 SO(3)  / \Omega S^3) \oplus \pi_k(S^2 \times S^2).$$
\end{thm}

\proof{}  First we define explicitly the map $h:P/\widetilde{P} \rightarrow \Im \Emb$.  
Let $\bar{a} \in P/\widetilde{P}$, i.e. $\bar{a}=a\widetilde{P}$, where $a \in P$.  We define $h(\bar{a})$ as the image of the composition $f(a) \circ i_c$  where $i_c$ is the standard symplectic embedding of $B^4(c)$ in $M_\mu$. It is obvious that $h$ is well defined because $f$ is an $H$--map and $f_{|\widetilde{P}}= \tilde{f}$ up to homotopy. Indeed, if we choose another representative of the same class $\bar{a} \in P/\widetilde{P}$, for example,  $\overline{ab}=ab\widetilde{P}$, with $b \in \widetilde{P}$, then   $h(\overline{ab})$ is given by the image of the composition $f(ab) \circ i_c$. Because $f$ is an H-map this composition is equal to $f(a) \circ f(b) \circ i_c \simeq f(a) \circ \widetilde{f}(b) \circ i_c$ up to homotopy. But $\tilde{f}(b) \in \Symp({\widetilde{M}}_{\mu,c})$ and this group is homotopy
 equivalent to the group of symplectic diffeomorphisms that preserve (not necessarily pointwise)  the image  $B_c$ of $i_c$, $\Symp(M_\mu,B_c)$ (see \cite[Lemmas 2.3 and 2.4]{LP}). It follows that up to homotopy $ \tilde{f}(b) \in \Symp(M_\mu,B_c)$. Hence the image of the composition $ f(a) \circ \tilde{f}(b) \circ i_c$ is the same as the image of $f(a) \circ i_c$ and we conclude that the map $h$ is well defined up to homotopy.

\smallskip

Therefore we obtain  the following homotopy commutative diagram 
$$\xymatrix{
\widetilde{P} \ \ \ar@{^{(}->}[r] \ar[d]^{\tilde{f}} &  P \ar[d]^f \ar[r] &  P/\widetilde{P}\ar[d]^h \\
\Symp({\widetilde{M}}_{\mu,c})  \ \ \ar@{^{(}->}[r] &  \Symp(M_\mu)\ar[r] & \Im \Emb
} $$
By the 5-lemma, since the first two vertical arrows are weak homotopy equivalences,  so is the last one. 

\medskip
 It remains to describe the quotient  $P/\widetilde{P}$ in more detail. Let $X \subset SO(3) \times SO(3)$ denote the image of a  section (we require that this section be a homomorphism, so $X$ is a subgroup of $SO(3) \times SO(3)$) of the principal fiber bundle 
$$ SO(3)\times SO(3)  \longrightarrow SO(3)$$ 
defined by quotienting $ SO(3)\times SO(3)$ by the diagonal $SO(3)$-action (for instance,  $\{\id \} \times SO(3)$ is such a section of $SO(3) \times SO(3)$). The base is not canonically identified with $SO(3)$ but, because the set of principal $SO(3)$--bundles over $SO(3)$ is trivial, we can find an identification making $\{\id \} \times SO(3)$ a section. Consider the following commutative diagram:
$$\xymatrix{
SO(3) \ar[r]^-{\Delta} \ar[d]^i & SO(3)\times SO(3) \ar[d]^{g_1} \\
SO(3) \times S^1 \ar[r]^-{g_2} &  SO(3) \times SO(3) \times S^1 }$$
where $g_1(a,b)=(a,b,1)$ and $g_2(c,d)=(c,c,d)$. The universal property of the pushout gives us a canonical continuous homomorphism 
\begin{equation}
\label{eq:ch}
P \longrightarrow SO(3)\times SO(3) \times S^1
\end{equation}
which is explicitely defined by mapping a sequence $x_1 \ldots x_n$ to the usual product in $SO(3)\times SO(3) \times S^1$ of the elements $x'_1 \ldots x'_n$ where if $x_i$ is equal to $(a,b)$ in $SO(3) \times SO(3)$, then $x'_i$ is $(a,b,1)$ and if $x_i$ is equal to $(c,t)$ in $SO(3) \times S^1$, then $x'_i$ is $(c,c,t)$.

Let $[S^1,X]$ denote the set of commutators $[c,x]=c^{-1}x^{-1}cx$ in $P$ with $c \in S^1-\{1\}$ and $x \in X-\{1\}$. Using the description we gave of the pushout of groups as a set of finite sequences with a certain equivalence relation it is not hard to see  that the set of commutators $[S^1,X]$ (this set is independent of the choice of section) is contained in the kernel of the canonical group homorphism (\ref{eq:ch}). Moreover, the kernel is the  free group generated by this set of commutators  (see a similar proof in \cite[Chapter 1,\S \,1]{Se}).  Thus we have a short exact sequence of topological groups 
$$ F[S^1,X] \longrightarrow P \longrightarrow SO(3)\times SO(3) \times S^1.$$

Similarly let $\tilde{X} \subset S^1 \times S^1$  denote the image of a section of the principal fiber bundle 
$$ S^1 \times S^1 \longrightarrow S^1$$ 
determined by quotienting $S^1 \times S^1$  by the diagonal $S^1$ action and consider the canonical group homomorphism 
$$ \widetilde{P} \longrightarrow S^1 \times  S^1 \times S^1.$$ We obtain another  short exact sequence of topological groups, in this case given by 
$$ F[S^1,\tilde{X}] \longrightarrow \widetilde{P} \longrightarrow S^1 \times S^1 \times S^1.$$

We can therefore write the following commutative diagram 
\begin{equation}
\label{seq}
\xymatrix{
F[S^1,\tilde{X}] \ar[r]^-{\tilde{\imath}} \ar@{^{(}->}[d]^{i_F} & \widetilde{P} \ar[r]^-{\tilde{\pi}} \ar@{^{(}->}[d]^{i_P} & S^1 \times  S^1 \times S^1 \ar@{^{(}->}[d]^{i_B}\\ F[S^1,X] \ar[r]^-i \ar[d]^{\pi_F} & P  \ar[r]^-\pi \ar[d]^{\pi_P} & SO(3) \times SO(3) \times S^1 \ar[d]^{\pi_B} \\ \frac{F[S^1,X]}{F[S^1,\tilde{X}]}  \ar[r]^-{\bar{\imath}} & P/\widetilde{P} \ar[r]^-{\bar{\pi}} & \frac{ SO(3) \times SO(3) \times S^1}{S^1 \times  S^1 \times S^1}}
\end{equation}

   It is easy to check that the two upper squares commute: that is to say $i_P \circ \tilde{\imath} = i \circ i_F$ and $i_B \circ \tilde{\pi} = \pi \circ i_P$.
The quotient spaces in the bottom row are Hausdorff spaces and the maps $\pi_B, \pi_F$ and $\pi_P$ are open and continuous. The maps $\bar{\pi}$ and $ \bar{\imath}$ are induced by the maps $\pi$ and $i$ respectively. If  $[\,a\,] \in P/\widetilde{P}$ and $[\, c\,] \in F[S^1,X]/F[S^1,\tilde{X}]$ then those  maps are defined by 
$$\bar{\pi} ([\, a\, ])= \bar{\pi}(a\widetilde{P})= (\pi_B \circ \pi)(a) ,$$
$$\bar{\imath}([\, c\, ])=\bar{\imath}(cF[S^1,\tilde{X}])=(\pi_P \circ i)(c).$$
They are well defined because  $\pi$ and  $i$ are group homomorphisms. Moreover, since the maps $i$, $\pi$, $\pi_F$, $\pi_P$ and $\pi_B$ are continuous, it follows easily that the maps $\bar{\pi}$ and $ \bar{\imath}$ are continuous too.

The bottom  row in  diagram (\ref{seq}) is the quotient of the two fibrations of groups on the first and second row. It is therefore a fibration too. We claim that  there exists a continuous section $\bar{\si}$ of $\bar{\pi}$. First note that both $\pi$ and $\tilde{\pi}$ have continuous sections $\si$ and $\tilde{\si}$ given by:
$$
\si(a,b,t) = (a,b) \,  (1,t) \qquad    \tilde{\si} (a,b,t) = (a,b) \,  (1,t).
$$

  This, by the way, implies that these two principal fibrations split topologically. However, the section $\si$ is not a homomorphism and therefore does not immediatly descend to a section of $\bar{\pi}$. But the restriction of $\si$ to the first two factors, i.e., to $SO(3) \times SO(3) \times \{1\}$ is obviously a group homomorphism, and therefore descends to a section
$$
\bar{\si}: \frac{ SO(3) \times SO(3)}{S^1 \times  S^1} = \frac{ SO(3) \times SO(3) \times S^1}{S^1 \times  S^1 \times S^1} \to P / \widetilde P.
$$
 Hence the long exact homotopy sequence of the fibration 
$$
\bar{\pi}:P/\widetilde P \to \frac{ SO(3) \times SO(3)}{S^1 \times  S^1}
$$
 breaks up into a family of short exact sequences
$$ 0 \to \pi_k\left(\frac{F[S^1,X]}{F[S^1,\tilde{X}]}\right) \to \pi_k(P/\widetilde P ) \to \pi_k \left(\frac{ SO(3) \times SO(3)}{S^1 \times  S^1}\right) \to 0.$$
Each of these sequence splits and therefore 
$$ \pi_k(P/\widetilde P ) \approx \pi_k\left(\frac{F[S^1,X]} {F[S^1,\tilde{X}]}\right) \oplus \pi_k \left(\frac{ SO(3) \times SO(3)}{S^1 \times  S^1}\right), \ \ \forall k \in \N.$$

\medskip
 In the quotient $(SO(3) \times SO(3))/(S^1 \times S^1)$ each $S^1$ factor is contained in the corresponding $SO(3)$ factor, hence the quotient is homeomorphic to the space $$SO(3)/S^1\times SO(3)/S^1 \approx S^2 \times S^2 .$$ 

Finally, it follows from the theory of amalgamated products of topological groups that the spaces  $[S^1,X] \subset P$ and $[S^1,\tilde{X}] \subset \widetilde{P}$ are  homeomorphic to $\Sigma SO(3)$ and $\Sigma S^1$ respectively. Writing $FY$ for the free topological group generated by the space $Y$, we have $FY \simeq \Omega \Sigma Y$. It follows that $F[S^1,X] /F[S^1,\tilde{X}] \simeq \Omega \Sigma^2 SO(3) /\Omega  S^3$. 

It remains to show  that the fibration 
\begin{equation} \label{mainfib}
\Omega \Sigma^2SO(3)/\Omega S^3 \rightarrow  P/\widetilde{P} \rightarrow S^2 \times S^2
\end{equation}
 is non trivial although the fibrations at the top and middle rows in diagram  (\ref{seq}) are. Note that  the quotient fibration above is not a principal bundle; the existence of a section implies the splitting of homotopy groups but not the splitting of the fibration. This non-splitting  is  an immediate consequence of the computation of the rational cohomology ring of  $\Im \Emb_{\om}(c, \mu)$ (cf.  Theorem \ref{thm:main} 1.b) or corollary \ref{cor:chring}).

\QED

From Theorem  \ref{th:main} it follows easily that:

\begin{thm} If $0<\mu-1 \leq c<1$ then the topological space $\Emb_{\om}(c, \mu)$ is weakly homotopy equivalent to the pull-back of the fibration 
$$ \Omega \Sigma^2 SO(3) /\Omega  S^3  \stackrel{\bar\imath}{\longrightarrow} P/\widetilde P \stackrel{\bar{\pi}}{\longrightarrow} S^2 \times S^2$$
by the fibered map $F_{\omega} \to S^2 \times S^2$ 
where $F_{\omega}$ is the space of symplectic frames over $M_\mu=S^2 \times S^2$.
\end{thm}

We leave the details of the proof of this result to the reader.

\section{The minimal models of $\Symp(M_{\mu}$)
 and of $\Symp(\tMuc)$ }\label{sc:mm1}

 First recall that in order to be applicable to some given topological space, the  theory of minimal models does not require that the space be simply connected. We just need that the space  has a nilpotent homotopy system which means that $\pi_1$ is nilpotent and $\pi_n$ is a nilpotent $\pi_1$--module for $n> 1$. Since the groups of symplectomorphisms $\Symp(M_{\mu}$) and  $\Symp(\tMuc)$ are $H$--spaces, it follows that they have a nilpotent homotopy system, because for a $H$--space $\pi_1$ is abelian and is therefore nilpotent, and moreover  $\pi_1$ acts trivially on all $\pi_n$'s. On the other hand, $\Im \Emb$ is simply connected because from computations in \cite{An,LP} we know that the generators of $\pi_1(\Symp(M_{\mu}))$ lift to the generators of $\pi_1(\Symp({\widetilde{M}}_{\mu,c}))$.  Therefore the theory of minimal models is applicable to all spaces under consideration.

   Recall that a model for a space $X$ is a graded differential algebra which provides a complete rational homotopy invariant of the space. Its cohomology is the rational cohomology of the space. The model can be constructed from the rational homotopy groups of $X$. In this case, it is always minimal, which implies that there is no linear term in the differential of the model, i.e the first term is quadratic. When there are no higher order term (i.e all terms are quadratic), then Sullivan's duality can be expressed in the following way:
$$
d b_k =  \sum_{i,j} \langle b_k,[b_i, b_j]\rangle    b_i b_j.
$$
where the $\langle a,b\rangle$ denotes the $a$-coefficient in the expression of $b$, and where the brackets denote the Whitehead product.  Finally, when $X$ is a H-space, as it is the case of both $\Symp(M_{\mu}$)
 and $\Symp(\tMuc)$, all Whitehead product vanish and the differential too.

   From these considerations and the computations of the rational homotopy groups of both $\Symp(M_{\mu}$) and of $\Symp(\tMuc)$ in \cite{Ab,LP}, we have:

\begin{itemize}

\item The minimal model of  $\Symp({\widetilde{M}}_{\mu,c})$  is $\La ({\tilde t},{\tilde x},{\tilde y}, {\tilde w})$, the free graded algebra generated by the elements ${\tilde t},{\tilde x},{\tilde y},{\tilde w}$ with degrees $\deg {\tilde t} = \deg {\tilde x} = \deg {\tilde y} = 1 $ and $ \deg {\tilde w} =2 $. 

\item The minimal model of $\Symp(M_{\mu})$ is $\La (t,x,y, w)$, the free graded algebra generated by the elements $t,x,y,w$ with degrees $\deg t = 1$, $\deg x = \deg y = 3$ and $ \deg w =4$.

\end{itemize}

   Let's now explain what these generators are. 

 The element $t$  in $\pi_1(\Symp(M_\mu))$ is the rotation in the fibers of the Hirzebruch surface $W_2$, once identified with $W_0 = S^2 \times S^2$; it is therefore the rotation in the fibers of $W_0 = B \times F$ round the two fixed symplectic surfaces in classes $B-F, B+F$ represented by the graph of the $\pm$ identity map from the base $B$ to the fiber $F$. The element $\tilde t$ is the blow-up of $t$ at the fixed point $([1,0], [0,0,1]) \in W_2$ identified with the center $\iota_c(0) \in S^2 \times S^2$ of the standard ball $B_c$. The element $x$ is the $3$-dimensional sphere generating $\pi_3(SO(3))$ where $SO(3)$ is considered as acting on the first factor in the obvious way, the element $y$ corresponds to the case when $SO(3)$ acts  on the second factor; the elements ${\tilde x},{\tilde y}$ are the blow-up of the $S^1$ part of that action that leaves the point  $\iota_c(0) \in S^2 \times S^2$ invariant. Finally, both $w$ and ${\tilde w}$ are symplectic elements that do not correspond to K\"ahlerian actions (i.e a symplectic action preserving  an integrable complex structure compatible with the symplectic form): $w$ is the Samelson product of $t$ and $x$, while ${\tilde w}$ is the Samelson product of $\tilde t$ and ${\tilde x}$.

\section{The minimal model of $\Im\Emb_{\om}(c, \mu)$}\label{sc:mm} \label{se:rathomtype}

   Any fibration $V \hookrightarrow P \to U$ for which the theory of minimal models applies (i.e. each space has a nilpotent homotopy system and the $\pi_1$ of the base acts trivially on the higher homotopy groups of the fiber) gives rise to a sequence
$$
( \La(U), d_U) \longrightarrow  (\La(U) \otimes \La(V), d) \longrightarrow (\La(V), d_V)
$$
where the differential algebra in the middle is a model for the total space. 
Let $d_{|U},d_{|V}$ represent the restriction of the differential $d$ to $U$ and $V$ respectively. The theory of minimal models implies that
$$ d_{|U} = d_U $$
$$ d_{|V} = d_V + d'$$
where $d'$ is a perturbation with image not in $\La (V)$

We have the fibration 
$$ \Symp({\widetilde{M}}_{\mu,c}) \longrightarrow \Symp(M_{\mu})  \longrightarrow \Im \Emb_{\om}(c,\mu) $$

\NI
and want to find the model for $\Im \Emb_{\om}(c,\mu)$. We know its algebra since  we know, from Lalonde-Pinsonnault \cite{LP}, its rational homotopy groups for $c \ge \mu -1$: 
$$ 
\pi_1 =0, \pi_2=\Q^2 ,  \pi_3= \Q^3 ,  \pi_4=\Q , \; \mbox{ and $\pi_n = 0$ for all $n>4.$}
$$ 
Therefore the algebra is $\La (a,b,e,f,g,h)$ where $ \deg e =\deg f = \deg g = 3$, $\deg a = \deg b = 2$  and $\deg h=4$.
Thus we get the following fibration
$$
\La (a,b,e,f,g,h), d_U \longrightarrow \La (a,b,e,f,g,h) \otimes \La (\tilde t,{\tilde x},{\tilde y}, {\tilde w}),d \longrightarrow \La (\tilde t,{\tilde x},{\tilde y}, {\tilde w}), d_V
$$
The differential $d$ satisfies $d_{|U} = d_U $ and $ d_{|V} = d_V+ d'=d'$. So in order to find the differential $d_U$ for $\La (a,b,e,f,g,h)$ it is sufficient to compute the differential $d$ for the model  $\La (a,b,e,f,g,h) \otimes \La (\tilde t,{\tilde x},{\tilde y}, {\tilde w})$. We need to compare this model with the minimal model  $\La (t,x,y, w)$ of $\Symp(M_{\mu})$ given in the last section.

\subsection{Computation of the differential d} 

Let us first apply the simplest method of dimension counting. That method yields easily the following partial results:

\begin{lemma}\label{le:diff} Without loss of generality, one may assume that the differential d satisfies:
$$
d \tilde  t = 0, \ \ d  {\tilde x} =a , \ \ d  {\tilde y}=b , \ \ d  {\tilde w} = g
$$
(and therefore the differentials of $a,b,g$ vanish). Moreover $d e$ and $d f$  must be quadratic, equal to (perhaps vanishing) linear combinations of $a^2, b^2, ab$. 
\end{lemma}

\proof{} Since the middle term computes the rational cohomology of $\Symp(M_{\mu})$, we need exactly one generator of degree  1. There is no loss of generality in assuming that it is $\tilde t$: $d \, \tilde  t=0$.
It follows that $d \, {\tilde x}$ and $d \, {\tilde y}$ must be different from 0 otherwise we would have too many generators in cohomology in dimension  1. By the theory of models for fibrations,  the perturbation $d'$ has image not  in $\La (\tilde t,{\tilde x},{\tilde y}, {\tilde w})$.
Therefore, without loss of generality, we may set:
$$d \, {\tilde x} =a , \ \ d \, {\tilde y}=b $$ which implies that $d \, a =d \, b =0$.

\noindent Now lets compute $d \, {\tilde w}$. It does not vanish because there is no generator in the cohomology of $\Symp(M_{\mu})$  in dimension 2. The theory of rational models for fibrations implies that the perturbation $d'$ is dual to the boundary operator $\partial: \pi_*(B) \otimes \Q \to \pi_*(F) \otimes \Q$. Since $\partial g = {\tilde w}$, we conclude that $d'{\tilde w} = g$, which means that $d{\tilde w}=g$ and implies that $dg=0$.

\noindent Now we need to determine the values of $d \, e $ and $d \, f$. They are linear combinations of$$0, \ a^2, \ b^2, \ ab, \ \mbox{and} \  h.$$
However, it is obvious that $h$ cannot appear in the differentials of $e,f$ since the coefficients of $h$ in both $d e $ and $d f$ would have to be the same in order to get a cycle $e-f$, but then one would only get a one-dimensional cohomology group in dimension $3$ instead  of a two-dimensional one.
\QED

\begin{Remark}  In what follows it might be useful to have in mind the following differentials:
$$ d(a \otimes {\tilde x})=a^2 , \  d(b \otimes {\tilde x})=ab , \  d(a \otimes {\tilde y})=ab , \ \mbox{and} \ d(b \otimes {\tilde y})=b^2 $$
and note that all choices of $a^2, b^2, ab$ for $d \, e$ and $d \, f$ lead to the same cohomology of $(\Lambda \otimes \Lambda, d + d')$, that is to say to the cohomology $\Lambda(t, x, y, w)$ with generators of dim 1,3,3,4 respectively. Indeed, one may quotient out $(\Lambda \otimes \Lambda, d + d')$ by the ideal generated by ${\tilde x},{\tilde y},a,b,{\tilde w},g,$ since all these elements kill each other.  Note that this ideal contains $a^2,b^2,ab$. Thus we get $\Lambda ({\tilde t}, e, f, h)$ with trivial differential, which is precisely what we are looking for, i.e it agrees with $\Lambda (t,x,y,w)$.  This shows that, using the ranks in cohomology alone, one cannot decide (among the choices $a^2,b^2,ab$) where $e,f$ are mapped to by the differential. And unfortunately, different choices lead to different cohomology for $\Lambda(a,b,e,f,g,h),d$. For instance, the choice $de = a^2, df = b^2$ is different from the choice
$de=a^2, df = b^2$ because in the second case, we get the cycle $e-f$ in dim 3. 
\end{Remark}

   To resolve the question, i.e to find the precise values of $d \, e, d \, f $ and $ d \, h$, lets first compute the Whitehead products $[a,a],$ $ [a,b], [b,b]$ in the rational homotopy of $\Im \Emb_{\om}(c,\mu)$. If the total space of the fibration $\Symp(\tMuc) \to \Symp(M_{\mu})  \to  \Im \Emb_{\om}(c, \mu)$ was contractible, computing such  products would boil down to computing the Samelson product of corresponding elements of the fiber. But our total space is not contractible, and we have to take also into account an horizontal part in the Whitehead product.
 
   Lets us briefly describe the generators of  $ \La(U)$, i.e. the generators of the rational homotopy groups of $\Im \Emb_{\om}(c,\mu)$. The group $\Symp(M_{\mu})$ acts on $\Im \Emb_{\om}(c,\mu)$ by $\phi \cdot A = \mbox{image} \, (\phi |_A)$ with stabilizer equal to $\Symp(M_{\mu},B_c)$, the subgroup of symplectic diffeomorphisms which preserve (not necessarily pointwise) $B_c$, the image of the standard embedding of the ball of capacity $c$ of $\R^{2n}$ in $M_{\mu}$. As we said in \S~\ref{se:intro}, this leads to the following homotopy fibration:
\begin{eqnarray}\label{fibration}
\Symp(\tMuc) \to \Symp(M_{\mu}) \to \Im \Emb_{\om}(c,\mu).  
\end{eqnarray}

The elements $e,f$ and $h$ are the images by the action of $\Symp(M_{\mu})$ on $\Im \Emb_{\om}(\mu,c)$ of the elements $x,y$ and $w$ of $\pi_*(\Symp(M_{\mu})) \otimes \Q$. 
The elements $a,b$  are uniquely defined as those spheres in the base of that fibration whose lifts to the total space $\Symp(M_{\mu}) $ are discs with boundary on the fiber equal to ${\tilde x}$ and ${\tilde y}$ respectively. These lifts are unique because $\pi_2(\Symp(M_{\mu})) \otimes \Q$ vanishes.
Finally, we choose the element $g$ in the following way. Its lift to the total space is a class in $ \pi_3(\Symp(M_{\mu}), \Symp(\tMuc)) \otimes \Q$  ($= \pi_3(\Im \Emb_{\om}(c, \mu)) \otimes \Q )$ whose boundary is equal to ${\tilde w}= [\tilde t,{\tilde x}]_S$ where $[ \cdot, \cdot]_S$ is the Samelson product.  Such a lift is not unique. To make it unique, we define it by first taking the $2$-disc $D_{\tilde x} \subset SO(3) = x \subset \Symp(M_{\mu})$ whose boundary is equal to $2 {\tilde x}$, and then taking the commutator of $t$ and $D_{\tilde x}$. This yields a $3$-disc $D$ lying inside $[t,x]_S = w$, whose boundary is the Samelson product $2[\tilde t,{\tilde x}] = 2{\tilde w}$. Set $g = D/2 \in \pi_3(\Symp(M_{\mu}), \Symp(\tMuc)) \otimes \Q$.

\begin{lemma} The Whitehead product $[a,b]$ vanishes, and 
$$
[a,a] = e \quad [b,b] = f.
$$
\end{lemma}

   \proof{} Consider the following commutative diagram:
$$
\begin{array}{ccc}
\Symp(\tMuc) &  & \Symp_p(M_{\mu})   \\
\downarrow &  & \downarrow  \\
\Symp(M_{\mu})  & = & \Symp(M_{\mu})  \\
\downarrow &  & \downarrow  \\
\Im \Emb_{\om}(c, \mu) & \stackrel{l}{\to} &  M
\end{array}
$$
where $\Symp_p(M_{\mu}) $ is the group of symplectic diffeomorphisms that fix a point $p$, the last horizontal map assigns to each unparametrised ball its center (well-defined up to homotopy),  and where the last downward arrow is the evaluation map at the point $p$. Write 
$$
[a,a] = \al e + \be f + \ga g.
$$
Applying $l_*$ to this equation yields
$$
[X_1, X_1] = \al Y_1 + \be Y_2  
$$
where $X,Y$ are the generators of $\pi_*(S^2) \otimes \Q$ and $X_i,Y_i$ their images in the $i^{th}$ factor of $S^2 \times S^2$. But $[X,X] = Y$ in the rational homotopy of $S^2$, so $\al = 1$ and $\be =0$. To compute $\ga$, recall that $a$ is defined as the element whose lift to  $\Symp(M_{\mu})$ is  ${\tilde x}$ and similarly, $g$ is the element whose lift is $\tilde w$. Hence the value of $\ga$ is given by:
$$
[{\tilde x},{\tilde x}]_S = \ga {\tilde w} .
$$ 
Because ${\tilde x}$ commutes with itself, $\ga$ vanishes. Hence $[a,a] = e$. The case of $[b,b]$ is similar, and the fact that $[a,b]$ vanishes is proved by the  same argument and the facts that $[X_1,X_2] = 0$ and  $[{\tilde x},{\tilde y}]_S = 0$. The first fact obviously holds because $[X_1,X_2]$ is the obstruction to extend the inclusion map $(S^2 \times \{ pt \}) \cup (\{ pt \} \times S^2) \to S^2 \times S^2$ to a map defined on $S^2 \times S^2$. The product $[{\tilde x},{\tilde y}]_S$  vanishes because ${\tilde x},{\tilde y}$ are $S^1$ actions of the same 2-torus, and therefore commute.  \QED

Recall that Sullivan's duality implies:
$$
d b_k =  \sum_{i,j} <b_k,[b_i, b_j]>   b_i b_j.
$$
Therefore, the last lemma implies that:
$$
d_U(e) = a^2  \quad  \mbox{and} \quad d_U(f) = b^2.
$$

 It remains to compute $dh$. 

\begin{lemma}\label{diffd}
The differential d satisfies $dh=k bg$ where $k$ is a non zero rational number. 
\end{lemma}
\proof{}
Notice that $dh \neq 0$ if and only if $H^4(\Im \Emb_{\om}(c, \mu); \Q)$ is one dimensional. Indeed, if $dh \neq 0$ there is only one element remaining in degree 4, namely $ab$. Hence, in that case,  we have $dh=g\tau$ where   $\tau$ is a non-zero linear combination of $a$ and $b$, since there are no closed classes  in degree 5 except for $g\tau$. Moreover, there is a constant $k\neq 0$ such that $dh=k bg$, because the Whitehead product $[a,g]$ vanishes. Indeed, recall the $2a$ is the projection on the base $\Im \Emb_{\om}(c, \mu)$ of the $2$--disc $D_{\tilde x} \subset D $ defined above, while $2g$ is the projection of the $3$--disc $D\subset [t,x]_S=w$. Therefore $a \subset g $ and their Whitehead product must vanish since $\pi_4(S^3) \otimes \Q=0$.    

Let $G$ denote the symplectomorphism group of $M_\mu$, $\Symp (M_\mu)$, and $K$ a suitable subgroup that can be identified with the group of symplectomorphisms of $M_\mu$ blown up at an embedded ball, $ \Symp(\tMuc)$. We will show that $H^4(\Im \Emb_{\om}(c, \mu); \Q)= H^4(G/K; \Q)$ is one dimensional by using an argument that uses the Eilenberg--Moore spectral sequence for the fibration $G/K \to BK \to BG$. This spectral sequence, which is a second quadrant spectral sequence, converges to $ H^*(G/K; \Q)$ and the $E_2$--term is given by 
$$ E_2^{i,j}=\Tor_{H^*(BG)}^{-i,j} (\Q\, , H^*(BK)).$$ 
We follow Paul Baum's paper \cite[Section 2]{Ba} to calculate these Tor groups. Let $\Lambda$ be a graded $\Q$--algebra and,  $M$ and $N$ be $\Lambda$--modules. Then  $\Tor_\Lambda(M,N)$ is the bigraded $\Q$--module obtained as follows. Consider  a projective resolution $R$ of $M$ over $\Lambda$ given by 
$$ R = \{ \xymatrix@+.5pc{ \hdots \ \ar[r] &  \ R^{(-2)} \ \ar[r]^{f^{(-2)}} &  \ R^{(-1)} \ \ar[r]^{f^{(-1)}} & \ R^{(0)} \ \ar[r]^{f^{(0)}} & \ M \ \ar[r] & \ 0 } \}. $$
Let $L$ be the bigraded differential $\Q$--module defined by $L^{p,q}=(R^{(p)} \otimes_\Lambda N )^q$ with $d: L^{p,q} \rightarrow  L^{p+1,q} $ given by $f^{(p)}\otimes_\Lambda 1_N$. $\Tor_\Lambda(M,N)$ is the homology of $L$, that is $\Tor_\Lambda^{p,q}(M,N)=H^{p,q}(L).$

In our example we have $ \Lambda=H^*(BG;\Q)$, $M=\Q$ and $N=H^*(BK;\Q)$. In order to calculate the cohomology rings of  $BG$ and $BK$  note  that  Theorem \ref{thm:htMuc} implies that the map $BK \rightarrow BG$ is equivalent to the map of homotopy pushouts of the following map of diagrams:
$$\xymatrix@C+2pc{ 
 &  BK_1 \ar[r]  & BG_1 \\
BK_0 \ar[ur] \ar[r] \ar[dr] & BG_0 \ar[ur] \ar[dr] & \\
& BK_2 \ar[r] & BG_2 }
$$
where $G_1=SO(3) \times SO(3)$, $G_2= SO(3) \times S^1$ and $G_0=SO(3)$ seen as the diagonal subgroup of $G_1$ and first factor in $G_2$. The groups $K_i \subset G_i$ are the standard maximal tori as  seen in Section \ref{se:fullhomtype}.
A simple calculation, using the Mayer--Vietoris sequence, then  gives the cohomology ring of the space $BK$ 
$$ H^*(BK;\Q)=\Q[z,r,s]/\langle z(r-s)\rangle, \mbox{ \ \ where \ \ } |r|=|s|=|z|=2,$$
and the cohomology of $BG$ can be identified with the subring 
$$ H^*(BG;\Q)=\Q[z,r^2,s^2]/\langle z(r^2-s^2)\rangle.$$
We need to construct a projective resolution for $\Q$ as a $H^*(BG)$--module. We can do that by means of the augmentation of $\Lambda$, $\varepsilon:\Q[z,r^2,s^2]/\langle z(r^2-s^2) \rangle \rightarrow \Q \simeq \Lambda^0 $. Therefore we may calculate these Tor groups using the following resolution (called the Koszul resolution)  
$$\Lambda (\alpha, \beta, \gamma, \delta) \otimes \Q[z,r^2,s^2]/\langle z(r^2-s^2) \rangle,$$
with differentials given by 
\begin{eqnarray}\label{diff}
d(\alpha)=z, \ d(\beta)=r^2, \ d(\gamma)=s^2, \ d(\delta)= \alpha(r^2-s^2).
\end{eqnarray}
Here $\Lambda(\alpha, \beta, \gamma, \delta)$ denotes the free (bi)graded algebra on elements $\alpha, \beta, \gamma$ and $\delta$ in bidegrees $(-1,2), (-1,4), (-1,4)$ and $(-2,6)$ respectively. The above complex is a module over $\Q[z,r^2,s^2]/\langle z(r^2-s^2) \rangle$  which is graded in external degree zero, i.e. lies in grading $(0,*)$. It follows that the Tor groups of interest are the cohomology of the complex 
$$\Lambda (\alpha, \beta, \gamma, \delta) \otimes \Q[z,r,s]/\langle z(r-s) \rangle.$$
Here we use the identification of $H^*(BG;\Q)$ as a subring of $ H^*(BK;\Q)$ and under this identification the differential of the complex above 
  satisfies the equalities (\ref{diff}) and $d(\eta \otimes u)= d\eta \otimes u $ with $\eta \in \Lambda (\alpha, \beta, \gamma, \delta)$ and  $u\in \Q[z,r,s]/\langle z(r-s)\rangle $. Any class in total degree 4, which is in negative external degree may be written as $x_4=c_1\delta+c_2\alpha\beta+c_3\alpha\gamma$, where $c_1,c_2$ and $c_3$ are constants. For it to be closed we need 
$$c_1\alpha(r^2-s^2)+z(c_2\beta+c_3\gamma)-\alpha(c_2r^2+c_3s^2)=0$$
which can happen only if $c_1=c_2=c_3=0$. Hence the only closed classes are in external degree zero, and therefore, the only class in Tor in total degree 4, is generated by $rs$.
\QED

 \bigskip

   This gives: 
\begin{thm}  \label{mmIEmb}The minimal model of $\Im \Emb_{\om}(c, \mu)$ is 
$$\La(\Im \Emb_{\om}(c, \mu))=\La (a,b,e,f,g,h)=\Lambda(S^2 \times S^2) \otimes \Lambda(g,h)$$ with generators in degrees $2,2,3,3,3,4$ and with differential
$$
d_U e = a^2, \; d_U f = b^2, \;  d_U g = d_U a = d_U b = 0, \; d_U h=kbg,
$$ 
where $\Lambda(S^2 \times S^2)$ is the minimal model for $S^2 \times S^2$ and $k$ is a non zero rational number. 
Then the rational cohomology ring of  $\Im \Emb_{\om}(c, \mu)$ is equal to the algebra 
$$ H^*(\Im \Emb_{\om}(c, \mu);\Q)=\Lambda(a,b,c,gh, \hdots, gh^n, \hdots,bh, \hdots, bh^n, \hdots)/\langle a^2,b^2,bg \rangle $$
where $n\in \N$ (see the computation of this cohomology ring in corollary \ref{cor:chring}).
 It is therefore not homotopy equivalent to a finite-dimensional CW-complex. 
\end{thm}

\section{The minimal model of $\Emb_{\om}(c,\mu)$} \label{se:rathomtypeb}

  In this section, we compute the minimal model of the space $\Emb_{\om}(c, \mu)$ of parametrised symplectic balls. Consider the fibration $U(2) \to \Emb_{\om}(c, \mu) \to \Im \Emb_{\om}(c, \mu)$.  First observe that this fibration is the restriction to $B_c$ of the fibration $\Symp(\tMuc) \to \Symp(M_{\mu}) \to \Im \Emb_{\om}(c, \mu)$. This can be expressed by the following commutative diagram:
$$
\begin{array}{ccccc}
\Symp^{\id, B_c}(M_{\mu})  &\hookrightarrow & \Symp^{U(2)}(M_{\mu}, B_c) & \stackrel{restr}{\to}  & U(2)   \\
| | & & \downarrow &  & \downarrow  \\
\Symp^{\id, B_c}(M_{\mu})  & \hookrightarrow & \Symp(M_{\mu})  & \stackrel{restr}{\to} &  \Emb_{\om}(c, \mu)  \\
\downarrow & &  \downarrow & & \downarrow  \\
\{B_c\} & \hookrightarrow & \Im \Emb_{\om}(c, \mu) & =  & \Im \Emb_{\om}(c, \mu) 
\end{array}
$$
where $restr$ is the restriction to the standard embedded ball $B_c \subset M_{\mu}$, $\Symp^{U(2)}(M_{\mu}, B_c)$ is the subgroup of $\Symp(M_{\mu})$ formed of diffeomorphisms that preserve the ball $B_c$ and act in a $U(2)$ linear way on it, and $\Symp^{\id, B_c}(M_{\mu})$ is the subgroup of $\Symp(M_{\mu})$ formed of the elements that fix the ball $B_c$ pointwise. Recall that there is a natural homotopy equivalence between $\Symp^{U(2)}(M_{\mu}, B_c)$  and $\Symp(\tMuc)$, so the vertical fibration in the middle is equivalent to the fibration $\Symp(\tMuc) \to \Symp(M_{\mu}) \to \Im \Emb_{\om}(c, \mu)$ in \cite{LP}. 

   We also have the commutative diagram:
$$
\begin{array}{ccc}
U(2) &  & U(2)   \\
\downarrow &  & \downarrow  \\
\Emb_{\om}(c, \mu)  &  \stackrel{j}{\to} & UFr(M)  \\
\downarrow &  & \downarrow  \\
\Im \Emb_{\om}(c, \mu) & \stackrel{l}{\to} &  M
\end{array}
$$
where $UFr(M)$ is the space of unitary frames of $M$, $j$ is the 1-jet map evaluated at the origin (followed by the Gram-Schmidt process assigning a unitary frame to each symplectic one), and where the last horizontal map assigns to each unparametrised ball its center (well-defined up to homotopy).  

    The minimal model for $U(2)$ is $\La(u_0,v_0)$ where $deg (u_0) = 1$ and $deg(v_0) = 3$. We first show that the elements $e,f,g,h \in \pi_*(\Im \Emb_{\om}(c, \mu)) \otimes \Q$ lift to $\pi_*(\Emb_{\om}(c, \mu)) \otimes \Q$, but not $a,b$. However the difference $a-b$ does lift. On the other hand only the element $v_0$ injects in $\pi_*(\Emb_{\om}(c, \mu)) \otimes \Q$, the element $u_0$ is killed.

\begin{prop} The rational homotopy of $\Emb_{\om}(c,\mu)$ is generated, as module over $\Q$, by a single element $\widetilde h$ in dimension $4$, by four elements $v,\widetilde e,  \widetilde f, \widetilde g$ in dimension $3$, and by a single element $\widetilde d_{a,b}$ in dimension $2$. The elements $\widetilde h, \widetilde e, \widetilde f$ are the images by the restriction map of the elements $w, x,y$ respectively. The element $v$ is the image of $v_0$, $\widetilde d_{a,b}$ is the unique lift of $a-b$ and $\widetilde g$ is a lift of $g$ well-defined up to a multiple of $v$. 
\end{prop}

\proof{} Consider the following commutative diagram of long exact sequences
$$
\begin{array}{ccccccccc}
\ldots  \to &  \pi_k(U(2)) \otimes \Q & \stackrel{\iota_*}{\to}  & \pi_k(\Emb_{\om}(c, \mu)) \otimes \Q  & \stackrel{\rho_*}{\to} & \pi_k(\Im \Emb_{\om}(c, \mu)) \otimes \Q & \stackrel{\partial_*}{\to}  \ldots  \\
  & \id \downarrow &  & j_*  \downarrow & & l_* \downarrow &    \\
 \ldots  \to &  \pi_k(U(2)) \otimes \Q & \to  & \pi_k(SFr(M)) \otimes \Q  & \to & \pi_k(M) \otimes \Q & \to  \ldots  
\end{array}
$$
Since $\pi_4(M) \otimes \Q$ vanishes, $l_{*=4}(h) = 0$, and therefore $\partial_{*}(h) =0$. Hence $\rho_{*=4}$ is an isomorphism between  $\pi_4(\Emb_{\om}(c, \mu)) \otimes \Q$ and $ \pi_4(\Im \Emb_{\om}(c, \mu)) \otimes \Q$. Let's denote by $\widetilde h$ the lift of $h$.

For $k=3$, the short sequence
$$
\pi_3(U(2)) \otimes \Q  \stackrel{\iota_*}{\to}   \pi_3(\Emb_{\om}(c, \mu)) \otimes \Q  \stackrel{\rho_*}{\to}  \pi_3(\Im \Emb_{\om}(c, \mu)) \otimes \Q 
$$
splits because, as we just saw $h$ is mapped to $0$, and $\pi_2(U(2))$ vanishes. Let's denote by $v$ the image of $v_0$ and by $\widetilde e, \widetilde f, \widetilde g$ the lifts of $e,f,g$,
 all are well-defined except $\widetilde g$ which is defined up to a multiple of the element $v$. Consider now the sequence
$$
 0 \to \pi_2(\Emb_{\om}(c, \mu)) \otimes \Q   \stackrel{\rho_*}{\to}  \pi_2(\Im \Emb_{\om}(c, \mu)) \otimes \Q \stackrel{\partial_*}{\to}  \pi_{1}(U(2)).
$$
The elements $a,b$ are by definition such that they lift to discs
$$
\phi_a, \phi_b: D^2 \to \Symp(M_{\mu})
$$
with boundary equal to the elements $x,y \in \pi_1(\Symp(\tMuc)) \otimes \Q$ respectively. Therefore, their lifts to $\Emb_{\om}(c, \mu) \otimes \Q$ are  the $2$-discs
$$
\psi_a, \psi_b: D^2 \to \Emb_{\om}(c, \mu)
$$ 
defined by $\psi_{a,b}(z) = \phi_{a,b} |_{B_c}$. Hence their boundaries are the restriction of the loops  ${\tilde x},{\tilde y} \in \pi_1(\Symp(M_{\mu},B_c)) \otimes \Q$ to the standard ball $B_c \subset M_{\mu}$. But each of these loops preserve $B_c$ (not pointwise) and correspond to the generator of $\pi_1(U(2)) \otimes \Q$ through the identification $B^4(c) (\subset \R^4) \to B_c$.
This proves that each of $a$ and $b$ is mapped to $u_0$ by the boundary operator of the above sequence. Denote by $\widetilde d_{a,b}$ the lift to $ \pi_2(\Emb_{\om}(c, \mu)) \otimes \Q$ of the element $d_{a,b}=a-b$.

   Finally, the map $\partial_*: \pi_2(\Im \Emb_{\om}(c, \mu)) \otimes \Q \to  \pi_{1}(U(2)) \otimes \Q$ being onto, the space $\pi_1(\Emb_{\om}(c, \mu)) \otimes \Q$ must vanish. \QED

\MS
  Let's compute the minimal model of $\Emb_{\om}(c, \mu)$.  By the last proposition, a model of $\Emb (c,\mu)$ is given by $(\La(\widetilde d_{a,b}, \widetilde e, \widetilde f, \widetilde g, v, \widetilde h), d_0)$. By minimality, there is no linear term in the differential, so $d_0(\widetilde d_{a,b}) = 0$, while the constants in the expression 
$$
d_0 (\widetilde e) =  c_1  \widetilde d_{a,b}^{\; 2}, \quad  d_0 (\widetilde f) =  c_2 \widetilde d_{a,b}^{\; 2} \quad d_0 (\widetilde g) =  c_3 \widetilde d_{a,b}^{\; 2} \quad d_0 v =  c_4 \widetilde d_{a,b}^{\; 2}
$$
are given, by duality, by
$$
[\widetilde d_{a,b}, \widetilde d_{a,b}] = c_1 \widetilde e + c_2 \widetilde f + c_3 \widetilde g + c_4 v.
$$ 
Denoting by $\rho$ the projection $\Emb_{\om}(c, \mu) \to \Im \Emb(c,\mu)$, we have:
$$
\rho_{*}([\widetilde d_{a,b},  \widetilde d_{a,b}]) = [d_{a,b}, d_{a,b}] = [ a-b, a-b] = [a,a] + [b,b] = e + f.
$$
Therefore $c_1 = c_2 = 1$ and $c_3=0$, and we get $[\widetilde d_{a,b}, \widetilde d_{a,b}] = \widetilde e +  \widetilde f  + c_4 v$. Now any value of this constant leads to the same model, up to isomorphism. Indeed, since $d_0 \widetilde e = d_0 \widetilde f = \widetilde d_{a,b}^{\; 2}$ and $d_0 v = c_4 \widetilde d_{a,b}^{\; 2}$, this means that $\widetilde e$ kills $\widetilde d_{a,b}^{\; 2}$ and thus both $\widetilde f$ and $v$ can be considered as cycles (up to a reparametrization of the basis of the algebra). 

   Finally, the differential of $\widetilde h$ is given by 
the coefficient affecting the term $\widetilde h$ in the Whitehead products $[\widetilde d_{a,b}, \widetilde e],
[\widetilde d_{a,b}, \widetilde f], [\widetilde d_{a,b}, \widetilde g], [\widetilde d_{a,b}, v]$. Projecting on the base of the fibration, we see that all these coefficients must vanish, except for the coefficient $k \in \Q$ in $d_0 \widetilde h = k \widetilde d_{a,b}\widetilde g$. Indeed projecting $[\widetilde d_{a,b}, \widetilde g]$ on the base we have 
$$
\rho_{*}([\widetilde d_{a,b},  \widetilde g]) = [ d_{a,b}, g]= [a-b,g]=[a,c]-[b,c]=-kh,
$$
since  $[a,g]=0$ and the differential $d$  of the minimal model of $\Im\Emb$ satisfies $dh=kbg$ for some $k \neq 0$, as seen  in  Lemma \ref{diffd}. This shows that the differential of $\widetilde h$ is given by 
$$ 
d_0\widetilde h=-k\widetilde d_{a,b} \widetilde g.
$$

  Denoting by $\widetilde f'$  and $v'$ the elements $\widetilde f - \widetilde e$ and $v -c_4 \widetilde e$ respectively, the set $\{\widetilde e, \widetilde f', \widetilde g, v'\}$ is a basis of the $3$-dimensional generators. This proves the following:

\begin{thm}\label{mmEmb}
A minimal model of $\Emb_{\om}{(c, \mu)}$ is $(\La(\widetilde d_{a,b}, \widetilde e, \widetilde f', \widetilde g, v', \widetilde h), d_0)$ with generators of degrees $2,3,3,3,3,4$ and with differential given by 
$$ d_0 \widetilde d_{a,b}  = d_0 \widetilde f' = d_0 \widetilde g = d_0 v'  = 0, \  d_0 \widetilde e =   \widetilde d_{a,b}^{\; 2}
 \mbox{ \  and \ } d_0\widetilde h=-k\widetilde d_{a,b} \widetilde g $$ where $k$ is a non zero rational number. Then the rational cohomology ring of $\Emb_{\om}{(c, \mu)}$ is given by 
$$ H^*(\Emb_{\om}{(c, \mu)}; \Q)=\Lambda (\widetilde d_{a,b},\widetilde f',\widetilde g, v', \widetilde g \widetilde h, \hdots,\widetilde g \widetilde h^n, \hdots, s_1, \hdots, s_n, \hdots)/\langle \widetilde d_{a,b}^{\;2},\widetilde d_{a,b}\widetilde g \rangle $$ where $s_n=\widetilde h^{n-1}(\widetilde h \widetilde d_{a,b} + nk \widetilde e \widetilde g)$ and $n \in \N$  (see the computation of this cohomology ring in corollary \ref{cohoemb}).
\end{thm}
\BS

\section{Cohomology rings } \label{se:cohomology}
It is easy  to describe the cohomology ring of $\Im \Emb_\omega(c,\mu)$ with rational coefficients.  A careful comparation between the Serre spectral sequence of the fibration 
\begin{equation}\label{eq:mf}
\Omega \Sigma^2 SO(3) /\Omega  S^3 \stackrel{\bar{\imath}}{\longrightarrow} P/{\widetilde{P}} \stackrel{\bar{\pi}}{\longrightarrow} S^2 \times S^2
\end{equation}
and Theorem \ref{mmIEmb}, together with Theorem \ref{th:main}, gives the  cohomology ring of $\Im \Emb_\omega(c,\mu)$.
\begin{cor}\label{cor:chring}
If $0<\mu-1\leq c< 1 $ the cohomology ring of $\Im \Emb_\omega(c,\mu)$ with rational coefficients is given by 
$$ H^*(\Im \Emb_\omega(c,\mu); \Q)= \Lambda(a,b,c,gh, \hdots, gh^n, \hdots,bh, \hdots, bh^n, \hdots)/ \langle a^2,b^2,bg \rangle, $$ 
that is, 
$$ H^*(\Im \Emb_\omega(c,\mu); \Q)= H^*(S^2 \times S^2;\Q) \otimes \Lambda(g,gh, \hdots, gh^n, \hdots,bh, \hdots, bh^n, \hdots)/ \langle bg \rangle, $$
where $n \in \N$, $f$ is a generator of $H^2(S^2 \times S^2;\Q)$, and $g,h$ correspond to the generators of the cohomology ring  $H^*( \Omega \Sigma^2 SO(3) /\Omega  S^3; \Q)$ where  $|g|=3$ and $|h|=4$.
\end{cor}
\proof{}
 First notice that at all odd primes, $SO(3) = S^3$ (more precisely there is only 2--torsion in $\pi_*(SO(3))$, so the localization, away from 2, of the two spaces is the same), hence the map $S^3 \to \Sigma^2 SO(3)$ is null homotopic away from the prime 2, and consequently that the fiber $\Omega \Sigma^2 SO(3) /\Omega  S^3$ is equivalent to the space $S^3 \times \Omega S^5$  away from the prime 2. Hence the cohomology ring of the fiber is given by $$H^*( \Omega \Sigma^2 SO(3) /\Omega  S^3; \Q)= \Lambda(g) \otimes \Q[h].$$ 
We showed in the proof of Lemma \ref{diffd} that $H^4(\Im \Emb_\omega(c,\mu),\Q)$ is one dimensional. This implies that in the $E_2$--term of the Serre spectral sequence of the fibration \eqref{eq:mf} the differential $d_2h$ does not vanish. More precisely, if we had $d_2h=0$ then $h$ would survive to the $E_\infty$ page of the spectral sequence and, unless $d_4g=ab$, we would have two generators in the cohomology group $H^4(\Im \Emb_\omega(c,\mu); \Q)$, namely $h$ and $ab$. However, we know that $d_rg=0$ for all $r \geq 2$ since $\rk H^3(\Im \Emb_\omega(c,\mu); \Q)=1$ as it is easy to verify from the minimal model. Therefore we can assume that  $d_2h=bg$.  It follows, for dimensional reasons, that the generators $bh^n$ where $n \in \N$ survive to the $E_\infty$ page  of the spectral sequence. A simple calculation then shows that the elements $gh^n$ where $n \in \N$ also survive to the $E_\infty$ page of the spectral sequence. Indeed the only way to kill these generators would be in the $E_4$ page if $d_4gh^n=abh^n$. However the computation of the minimal model of $\Im \Emb_\omega(c,\mu)$ implies  that  $\rk H^{4n+3}(\Im \Emb_\omega(c,\mu); \Q)=1$, and it is easy to verify that $gh^n$, for each $n$,  is the single element in dimension $4n+3$ that can survive to the $E_\infty$--page of the spectral sequence. Therefore $d_4gh^n$ vanishes for all $n \in \N$ and this completes the proof. \QED

 Using a similar argument we can compute the cohomology ring $\Im \Emb$ with $\Z_p$  coefficients with $p$ prime. Let $\Gamma_{\Z_p}[x]$ denote  the divided polynomial algebra on the generator $x$. This is defined to be the free $\Z_p$--module with basis $x_0=1,x_1,x_2,\hdots$ and multiplication defined by $x_ix_j= \left(\begin{array}{c} i+j \\ i \end{array} \right) x^{i+j}$. Moreover, there is an isomorphism $\Gamma_{\Z_p}[x] \approx \Z_p[x_1,x_p,x_{p^2},\hdots]/\langle x_1^p,x_p^p,x_{p^2}^p,\hdots\rangle =\bigotimes_{i\geq 0}\Z_p[x_{p^i}]/\langle x_{p^i}^p\rangle $.
\begin{cor} \label{co:crzp}If $0<\mu-1\leq c< 1 $ and $p \neq 2$ then 
$$ H^*(\Im \Emb_\omega(c,\mu); \Z_p)= \Lambda(a,b,g)/\langle a^2,b^2,bg \rangle \otimes g\Gamma_{\Z_p}[h] \otimes b\Gamma_{\Z_p}[h]$$
where $|a |= |b |= 2$, $|g|=3$, $|h|=4$ and $\tau\Gamma_{\Z_p}[h]$, with $\tau = g$ or $\tau = b$, stands for the infinitely generated algebra $\Z_p[\tau h_1,\hdots, \tau h_1^{p-1},\tau h_p,\hdots, \tau h_p^{p-1},\tau h_{p^2}, \hdots,\tau h_{p^2}^{p-1}, \hdots ]$ where the generators $h_i$ are the generators of the divided polynomial algebra $\Gamma_{\Z_p}[h]$ as described above.
\end{cor}
\proof{}
As noticed in the proof of the previous corollary the fiber $\Omega \Sigma^2 SO(3) /\Omega  S^3$ is equivalent to the space $S^3 \times \Omega S^5$  away from the prime 2. Therefore we get that $$H^*( \Omega \Sigma^2 SO(3) /\Omega  S^3; \Z_p)= \Lambda(g) \otimes \Gamma_{\Z_p}[h],$$ where $p\neq 2$, $|g| = 3$ and $|h| = 4$. The same argument of the proof of the Lemma \ref{diffd}, using the Eilenberg--Moore spectral sequence, shows that $H^4(\Emb_\omega(c,\mu); \Z_p)$ is one dimensional if $p\neq 2$.  Since $\rk H^{4n+3}(\Im \Emb_\omega(c,\mu); \Q)=1$ where $n \in \N_0$ it follows that $ H^{4n+3}(\Im \Emb_\omega(c,\mu); \Z_p)$ is at least one dimensional. Then using again the Serre spectral sequence of fibration \eqref{eq:mf} and an argument similar to the one used in the above corollary we obtain the desired result. \QED

Next we will see that $\Im \Emb$ has $\Z_2$--torsion and therefore  the cohomology ring with these  coefficients is not as simple to describe as  the previous ones.
\begin{cor} \label{co:crz2}  When   $0<\mu-1\leq c< 1 $, the homology groups with $\Z_2$ coefficients of the space $\Im \Emb_\omega(c,\mu)$  are given by 
\begin{equation}\label{isoh} H_*(\Im \Emb_\omega(c,\mu); \Z_2) =  H_*(S^2\times S^2; \Z_2) \otimes  H_*(\Omega \Sigma^2 SO(3)/\Omega S^3;\Z_2).
\end{equation}
Moreover, as an algebra, we have that  $$H_*(\Omega \Sigma^2 SO(3)/\Omega S^3;\Z_2)=T(w_2,w_3,w_4) \otimes_{T(w_2)}\Z_2$$ where $T$ denotes the tensor algebra, that is, the free noncommutative algebra on the generators $w_i$ and  $|w_i|=i$. Therefore 
the cohomology ring of $\Im \Emb_\omega(c,\mu)$ with $\Z_2$ coefficients is given by 
$$ H^*(\Im \Emb_\omega(c,\mu); \Z_2) \cong   H^*(S^2\times S^2; \Z_2) \otimes  A$$ where $A$ has an infinite number of generators.
\end{cor}
\proof{}
Since the inclusions $i_F, i_P$ and $i_B$ of diagram (\ref{seq}) induce injective maps in homology with $\Z_2$ coefficients, it follows from the Leray--Hirsch Theorem that we have the following isomorphisms as vector spaces
$$H^*( \Omega \Sigma^2 SO(3);\Z_2) \cong H^*( \Omega \Sigma^2 SO(3) /\Omega  S^3; \Z_2) \otimes H^*(\Omega  S^3; \Z_2), $$
$$ H^*(P ; \Z_2) \cong H^*(P/{\widetilde{P}} ; \Z_2)\otimes H^*({\widetilde{P}} ; \Z_2) \ \ \ {\mbox{and}} $$
$$H^*(S^1 \times SO(3) \times SO(3); \Z_2) \cong H^*(S^2\times S^2; \Z_2) \otimes H^*(S^1 \times S^1 \times S^1; \Z_2).$$
Moreover, since the fibrations $\pi:P \rightarrow S^1 \times SO(3) \times SO(3)$ and $\tilde{\pi}:\widetilde{P} \rightarrow S^1 \times S^1 \times S^1$ are (weak) homotopically trivial then we obtain the following isomorphisms as graded algebras 
$$  H^*(P ; \Z_2) \cong H^*(S^1 \times SO(3) \times SO(3); \Z_2)\otimes H^*( \Omega \Sigma^2 SO(3);\Z_2), $$
$$H^*({\widetilde{P}} ; \Z_2) \cong H^*(S^1 \times S^1 \times S^1; \Z_2) \otimes H^*(\Omega  S^3; \Z_2). $$
The 5 previous isomorphisms yield the isomorphisms  
\begin{equation}\label{ve:iso2}
H^*(P/{\widetilde{P}} ; \Z_2) \cong H^*(\Omega \Sigma^2 SO(3)/\Omega S^3;\Z_2) \otimes  H^*(S^2\times S^2; \Z_2)
\end{equation}
and \eqref{isoh} as vector spaces. 

It follows that the homomorphism   ${\bar{\imath}}^*$ is surjective  and the Serre spectral sequence of the fibration (\ref{eq:mf}) collapses at $E_2$.
Therefore
$$E_\infty^{*,*} \cong E_2^{*,*} \cong H^*(\Omega \Sigma^2 SO(3)/\Omega S^3;\Z_2) \otimes  H^*(S^2\times S^2; \Z_2)$$ 
as bigraded algebras. But this does not directly shows that the isomorphism (\ref{ve:iso2})  also holds as a graded algebra isomorphism. However, it is clear that $H^*(P/{\widetilde{P}} ; \Z_2)$  has a subalgebra ${\bar{\pi}}^*(H^*(S^2\times S^2; \Z_2)) \cong H^*(S^2\times S^2; \Z_2)$. 
Although it is not easy to describe the $\Z_2$--cohomology of the space $\Omega \Sigma^2 SO(3) /\Omega  S^3$  one can calculate its $\Z_2$--homology. To to this recall that the homology of $\Omega \Sigma X$ is a Tensor algebra on the homology of $X$ for any connected space $X$. Hence the map $\Omega  S^3 \to \Omega \Sigma^2 SO(3)$  corresponds to the obvious inclusion of tensor algebras over $\Z_2$:
$$ T(w_2) \to T(w_2,w_3,w_4)$$ where $T$ denotes the tensor algebra and  $|w_i|=i$. A spectral sequence known as the Bar spectral sequence for a principal fibration, can then be applied to give 
$$ H_*( \Omega \Sigma^2 SO(3) /\Omega  S^3; \Z_2)= T(w_2,w_3,w_4) \otimes_{T(w_2)}\Z_2.$$
In \cite{KLW} the reader will find the necessary results on the Bar spectral sequence (cf. Theorem 4.2) and further references. 
It follows that  the  cohomology ring  $H^*(\Omega \Sigma^2 SO(3) /\Omega  S^3; \Z_2)$ has an infinite number of generators. From the $E_2$ page of the spectral sequence of the fibration (\ref{eq:mf}) we can then conclude that 
$$H^*(P/{\widetilde{P}} ; \Z_2) \cong   H^*(S^2\times S^2; \Z_2) \otimes A, $$
 as graded algebras, where $A$ has an infinite number of generators, but it is not necessarily isomorphic as a graded  algebra to $H^*(\Omega \Sigma^2 SO(3)/\Omega S^3;\Z_2)$. This isomorphism together with Theorem \ref{th:main} completes the proof. 
\QED

A careful comparation of Theorem \ref{mmEmb} and the Serre spectral sequence of the fibration
\begin{equation}\label{fEmb}
 U(2) \longrightarrow  \Emb_\omega(c,\mu) \longrightarrow \Im \Emb_\omega(c,\mu)
\end{equation}  yields the cohomology ring of $\Emb_\omega(c,\mu)$ with rational coefficients.
\begin{cor} \label{cohoemb} If $0<\mu-1\leq c< 1 $ then 
$$H^*( \Emb_\omega(c,\mu); \Q)\cong \Lambda( b, f, g,v, gh, \hdots, gh^n, \hdots,bh, \hdots, bh^n, \hdots)/\langle b^2,bg \rangle$$ where $H^*(U(2);\Q) \cong \Lambda(u,v)$, $|b|=3$ and  $a,b,g,gh^n,bh^n$ with $n\in \N$ correspond to  the generators of the cohomology ring of $\Im \Emb_\omega(c,\mu)$. \end{cor}
\proof From the computation of the minimal model of the space $ \Emb_\omega(c,\mu) $ in Theorem \ref{mmEmb} it follows that there is no generator in degree 1 in its cohomology ring. Hence in the $E_2$--page of the Serre spectral sequence of the fibration \eqref{fEmb} the differential $d_2$ satisfies $d_2u \neq 0$.
 Therefore $d_2u$ is a linear combination of $a$ and $b$. Notice that  the minimal model computation also shows that there is no element in degree 4 in the cohomology ring so the element $ab$ in the $E_2$--page has to be in the image of $d_2$.  We can choose $d_2u=a$ and then we have $d_2ub=ab$ as desired. The computation of the minimal model implies  that the element $ua$ survives to the $E_\infty$ page, as well as the generators $v$ and $g$, since we should have three generators of degree 3 in the cohomology ring. The element $ua$ correponds to the generator $f$. It is not hard to see that the generators $gh^n$ also survive  to the  $E_\infty$ page and they correspond to the generators $\widetilde g \widetilde h^n$ in the minimal model. Finally we see that the generators $bh^n$ cannot be in the image of $d_r$ with $r \geq 2$, so they also survive to the $E_\infty$ page. Moreover they correspond to the elements $s_n=\widetilde h^{n-1}(\widetilde h \widetilde d_{a,b} + nk \widetilde e \widetilde g)$ in the minimal model, where $n \in \N$, which clearly satisfy $d_0 s_n=0$.  
\QED
\begin{rmk}
Notice that this cohomology ring is equivalent to the one given in Theorem \ref{mmEmb} and Theorem \ref{thm:main}. Indeed the difference between the two is that here we use the generators of the cohomology ring of $\Im \Emb_\omega(c,\mu)$ to describe the ring while  there we use the generators of the minimal model.
\end{rmk}
Since there is no $\Z_p$--torsion if $p\neq 2$ it follows that we have a similar result for the cohomology with $\Z_p$ coefficients when $p\neq 2$.
\begin{cor} If $0<\mu-1\leq c< 1 $ and $p\neq 2$ then 
$$H^*( \Emb_\omega(c,\mu); \Z_p)\cong \Lambda(b,f,g,v)/\langle b^2,bg \rangle \otimes  g \,\Gamma_{\Z_p}[h]\otimes  b\Gamma_{\Z_p}[h],$$
where $\tau\Gamma_{\Z_p}[h]$, with $\tau = g$ or $\tau = b$, stands for the infinitely generated algebra $$\Z_p[\tau h_1,\hdots, \tau h_1^{p-1},\tau h_p,\hdots, \tau h_p^{p-1},\tau h_{p^2}, \hdots,\tau h_{p^2}^{p-1}, \hdots ]$$ where the generators $h_i$ are the generators of the divided polynomial algebra $\Gamma_{\Z_p}[h]$.
 
\end{cor}
In the case of $\Z_2$ coefficients is again more difficult to describe the cohomology ring since there is  $\Z_2$--torsion in $\Im \Emb_\omega(c,\mu)$ as we showed in Corollary \ref{co:crz2}. Using this corollary and the Serre spectral sequence of fibration \eqref{fEmb} we obtain
\begin{cor}
If $0<\mu-1\leq c< 1 $ then 
$$H^*( \Emb_\omega(c,\mu); \Z_2)\cong \Lambda(b,f,v)/\langle b^2 \rangle  \otimes A'$$
where the algebra $A'$  has an infinite number of generators.
\end{cor}

\end{document}